\documentclass[a4, 12pt]{amsart}
\usepackage[mathscr]{eucal}
\usepackage{amssymb}
\usepackage{latexsym}
\usepackage{amsthm}
\theoremstyle{plain}
\newtheorem{theorem}{Theorem}[section]
\newtheorem{proposition}{Proposition}[section]
\newtheorem{corollary}{Corollary}[section]
\newtheorem{remark}{Remark}[section]
\newtheorem{lemma}{Lemma}[section]
\newtheorem{example}{Example}[section]
\newtheorem{definition}{Definition}[section]

\makeatletter
\@addtoreset{equation}{section}

\title[Rigidity theorems of $\lambda$-hypersurfaces]
{Rigidity theorems of $\lambda$-hypersurfaces}
\author{Qing-Ming Cheng, Shiho Ogata and  Guoxin Wei}
\address{Qing-Ming Cheng \\ Department of Applied Mathematics, Faculty of Sciences,
Fukuoka  University, 814-0180, Fukuoka,  Japan, cheng@fukuoka-u.ac.jp}
\address{Shiho Ogata \\ Department of Applied Mathematics, Graduate School of Sciences,
Fukuoka  University, 814-0180, Fukuoka,  Japan}
\address{Guoxin Wei \\  School of Mathematical Sciences, South China Normal University,
510631, Guangzhou,  China, weiguoxin@tsinghua.org.cn}

\begin{document}
\maketitle

\begin{abstract}
\noindent Since  $n$-dimensional $\lambda$-hypersurfaces in the Euclidean space $\mathbb {R}^{n+1}$
are  critical points of the weighted area functional for
the weighted volume-preserving variations, in this paper, we study the rigidity properties of  complete $\lambda$-hypersurfaces.
We give a gap theorem of complete $\lambda$-hypersurfaces with polynomial area growth. By making use of
the generalized maximum principle for $\mathcal L$  of  $\lambda$-hypersurfaces, we prove  a  rigidity theorem of
complete $\lambda$-hypersurfaces.
\end{abstract}



\footnotetext{The first author was partially  supported by JSPS Grant-in-Aid for Scientific Research (B): No. 24340013
and Challenging Exploratory Research No. 25610016.
The third author was partly supported by grant No. 11371150 of NSFC.}

\section {Introduction}

\noindent
Let $X: M\rightarrow \mathbb{R}^{n+1}$ be a smooth $n$-dimensional immersed hypersurface in the $(n+1)$-dimensional
Euclidean space $\mathbb{R}^{n+1}$. In \cite{cw2}, Cheng and Wei have introduced notation of
the weighted volume-preserving mean curvature flow, which is defined as the following:
a family $X(\cdot, t)$ of smooth immersions
$$
X(\cdot, t):M\to  \mathbb{R}^{n+1}
$$
with  $X(\cdot, 0)=X(\cdot)$ is called {\it a weighted volume-preserving mean curvature flow} if
\begin{equation}
\dfrac{\partial X(t)}{\partial t}=-\alpha(t) N(t) +\mathbf{H}(t)
\end{equation}
holds, where
$$
\alpha(t) =\dfrac{\int_MH(t)\langle N(t), N\rangle e^{-\frac{|X|^2}2}d\mu}{\int_M\langle N(t), N\rangle e^{-\frac{|X|^2}2}d\mu},
$$
 $\mathbf{H}(t)=\mathbf{H}(\cdot,t)$ and $N(t)$ denote the mean curvature vector  and the normal vector of hypersurface
 $M_t=X(M^n,t)$ at point $X(\cdot, t)$, respectively  and   $N$ is the unit normal vector of $X:M\to  \mathbb{R}^{n+1}$.
One  can prove  that the flow (1.1) preserves the weighted volume $V(t)$ defined by
$$
V(t)=\int_M\langle X(t),N\rangle e^{-\frac{|X|^2}{2}}d\mu.
$$
\noindent
{\it The weighted area functional}
$A:(-\varepsilon,\varepsilon)\rightarrow\mathbb{R}$ is defined by
$$
A(t)=\int_Me^{-\frac{|X(t)|^2}{2}}d\mu_t,
$$
where $d\mu_t$ is the area element of $M$ in the metric induced by $X(t)$.
Let $X(t):M\rightarrow \mathbb{R}^{n+1}$ with $X(0)=X$ be a  variation of $X$.
If $V(t)$ is constant for any $t$, we call  $X(t):M\rightarrow \mathbb{R}^{n+1}$ {\it  a weighted volume-preserving variation of $X$}.
Cheng and Wei \cite{cw2} have proved that  $X:M\rightarrow \mathbb{R}^{n+1}$ is a critical point of the weighted area functional $A(t)$
for all weighted volume-preserving variations if and only  if there exists constant $\lambda$ such that
\begin{equation}
\langle X, N\rangle +H=\lambda.
\end{equation}
An immersed hypersurface $X(t):M\rightarrow \mathbb{R}^{n+1}$ is called {\it a $\lambda$-hypersurface} if the equation (1.2) is satisfied.

\begin{remark}
If $\lambda=0$, then the $\lambda$-hypersurface is a self-shrinker of the mean curvature flow. Hence, the $\lambda$-hypersurface is
a generalization of the self-shrinker.
\end{remark}
\begin{example}
The $n$-dimensional sphere $S^n(r)$ with radius $r>0$ is a compact $\lambda$-hypersurface in $\mathbb{R}^{n+1}$
with $\lambda=\frac{n}r-r$.
\end{example}
\begin{example}
For $1\leq k\leq n-1$, the $n$-dimensional cylinder  $S^k(r)\times \mathbb{R}^{n-k}$ with radius $r>0$ is a complete
and non-compact $\lambda$-hypersurface in $\mathbb{R}^{n+1}$ with $\lambda=\frac{k}r-r$.
\end{example}
\begin{example}
The $n$-dimensional Euclidean space $\mathbb{R}^{n}$ is a complete and non-compact $\lambda$-hypersurface
in $\mathbb{R}^{n+1}$ with $\lambda=0$.
\end{example}
\begin{definition}
If $X: M\rightarrow \mathbb{R}^{n+1}$ is an $n$-dimensional hypersurface in $\mathbb{R}^{n+1}$,
we say that $M$ has polynomial area growth if there exist constant $C$ and $d$ such that for all $r\geq 1$,
\begin{equation}
{\rm Area}(B_r(0)\cap X(M))=\int_{B_r(0)\cap X(M)}d\mu\leq Cr^d,
\end{equation}
where $B_r(0)$ is a standard ball in $\mathbb{R}^{n+1}$ with radius $r$ and centered at the origin.
\end{definition}

\noindent
In \cite{cw2}, Cheng and Wei have  studied  properties of  complete $\lambda$-hypersurfaces  with polynomial area growth.
They have proved that  a complete and non-compact $\lambda$-hypersurface  $X: M\rightarrow \mathbb{R}^{n+1}$
in the Euclidean space $\mathbb{R}^{n+1}$ has polynomial area growth if and only if  $X: M\rightarrow \mathbb{R}^{n+1}$
is a complete proper hypersurface.  Furthermore, there is a positive constant $C$ such that for $r\geq1$,
\begin{equation}
{\rm Area}(B_r(0)\cap X(M))=\int_{B_r(0)\cap X(M)}d\mu\leq Cr^{n+\frac{\lambda^2}2-2\beta-\frac{\inf H^2}2},
\end{equation}
where  $\beta=\frac{1}{4}\inf(\lambda-H)^2$.

\vskip 5mm
\noindent
In this paper, we study the rigidity of complete $\lambda$-hypersurfaces. We will prove the following:

\begin{theorem}
Let $X:M\rightarrow \mathbb{R}^{n+1}$ be an $n$-dimensional complete $\lambda$-hypersurface
with polynomial area growth  in the Euclidean space $\mathbb{R}^{n+1}$. Then  $X:M\rightarrow \mathbb{R}^{n+1}$
satisfies one of the following:
\begin{enumerate}
\item  $X:M\rightarrow \mathbb{R}^{n+1}$ is isometric to the sphere $S^n(r)$ with radius $r>0$,
\item $X:M\rightarrow \mathbb{R}^{n+1}$ is isometric to the Euclidean space  $ \mathbb{R}^{n}$,
\item $X:M\rightarrow \mathbb{R}^{n+1}$ is isometric to the cylinder $S^1(r)\times   \mathbb{R}^{n-1}$,
\item $X:M\rightarrow \mathbb{R}^{n+1}$ is isometric to the cylinder $S^{n-1}(r)\times   \mathbb{R}$,
\item $X:M\rightarrow \mathbb{R}^{n+1}$ is isometric to the cylinder $S^{k}(\sqrt{k})\times   \mathbb{R}^{n-k}$ for $2\leq k\leq n-2$,
\item there exists $p\in M$ such that the squared norm $S$ of the second fundamental form and the mean
curvature $H$ of $X:M\rightarrow \mathbb{R}^{n+1}$ satisfy
\begin{equation}
\biggl(\sqrt {S(p)-\dfrac{H^2(p)}n}+|\lambda| \dfrac{n-2}{2\sqrt{n(n-1)}}\biggl)^2+ \frac1n(H(p)-\lambda)^2> 1+\dfrac{n\lambda^2}{4(n-1)}.
\end{equation}
\end{enumerate}
 \end{theorem}

\begin{corollary}
Let $X:M\rightarrow \mathbb{R}^{n+1}$ be an $n$-dimensional complete $\lambda$-hypersurface
with polynomial area growth  in the Euclidean space $\mathbb{R}^{n+1}$.
If the squared norm $S$ of the second fundamental form and the mean curvature $H$ of $X:M\rightarrow \mathbb{R}^{n+1}$ satisfies
\begin{equation}\label{eq:001}
\biggl(\sqrt {S-\dfrac{H^2}n}+|\lambda| \dfrac{n-2}{2\sqrt{n(n-1)}}\biggl)^2+ \frac1n(H-\lambda)^2\leq 1+\dfrac{n\lambda^2}{4(n-1)},
\end{equation}
then $X:M\rightarrow \mathbb{R}^{n+1}$ is isometric to one of the following:
\begin{enumerate}
\item the sphere $S^n(r)$ with radius $0<r\leq\sqrt{n}$,
\item the Euclidean space  $\mathbb{R}^{n}$,
\item the cylinder $S^1(r)\times \mathbb{R}^{n-1}$ with radius $r>0$ and $n=2$ or with radius $r\geq 1$ and $n>2$,
\item the cylinder $S^{n-1}(r)\times \mathbb{R}$ with radius $r>0$ and $n=2$ or with radius $r\geq\sqrt{n-1}$ and $n>2$,
\item the cylinder $S^{k}(\sqrt{k})\times   \mathbb{R}^{n-k}$ for $2\leq k\leq n-2$.   
\end{enumerate}
 \end{corollary}
\begin{remark} If $\lambda=0$, that is, $X:M\rightarrow \mathbb{R}^{n+1}$ is an $n$-dimensional complete self-shrinker,
our condition \eqref{eq:001} becomes $S\leq 1$. Hence, our theorem is a general  generalization of Cao and Li \cite{[CL]}  and Le and  Sesum \cite{[LS]}
to $\lambda$-hypersurfaces. On study of complete self-shrinkers,  see \cite{cp},  \cite{cw}, \cite{[CZ]}, \cite{[CM]},  \cite{[DX1], [DX2]}, \cite{[H2], [H3]}.
\end{remark}
\begin{theorem}
Let $X:M\rightarrow \mathbb{R}^{n+1}$ be an $n$-dimensional complete $\lambda$-hypersurface
with polynomial area growth  in the Euclidean space $\mathbb{R}^{n+1}$. If
\begin{equation}
\biggl(H-\frac{\lambda}{2}\biggl)^2\geq n+\frac{\lambda^2}{4},
\end{equation}
then $(H-\frac{\lambda}{2})^2\equiv n+\frac{\lambda^2}{4}$ and $M$ is isometric to the sphere $S^n(r)$ with radius $r>0$.
\end{theorem}

\noindent
If we do not assume that $X:M\rightarrow \mathbb{R}^{n+1}$ has polynomial area growth, we can prove
the following:
\begin{theorem}
Let $X:M\rightarrow \mathbb{R}^{n+1}$ be an $n$-dimensional complete $\lambda$-hypersurface
in the Euclidean space $\mathbb{R}^{n+1}$.
If the squared norm $S$ of the second fundamental form and the mean curvature $H$ of $X:M\rightarrow \mathbb{R}^{n+1}$ satisfy
\begin{equation}\label{eq:003}
\sup\bigl\{(\sqrt {S-\dfrac{H^2}n}+|\lambda| \dfrac{n-2}{2\sqrt{n(n-1)}})^2+\frac1n(H-\lambda)^2\bigl\}< 1+\dfrac{n\lambda^2}{4(n-1)},
\end{equation}
then $X:M\rightarrow \mathbb{R}^{n+1}$ is isometric to one of the following:
\begin{enumerate}
\item the sphere $S^n(r)$ with radius $r<\sqrt{n}$,
\item the Euclidean space  $ \mathbb{R}^{n}$.
\end{enumerate}
 \end{theorem}
\begin{remark} If $\lambda=0$, that is, $X:M\rightarrow \mathbb{R}^{n+1}$ is an $n$-dimensional complete self-shrinker,
our condition \eqref{eq:003} becomes $\sup S< 1$. Our theorem  is a general  generalization of Cheng and Peng \cite{cp}
to $\lambda$-hypersurfaces.
\end{remark}
\noindent
We next give the following:

\begin{proposition}
Let $X:M\rightarrow \mathbb{R}^{n+1}$ be an $n$-dimensional compact $\lambda$-hypersurface in the Euclidean space $\mathbb{R}^{n+1}$. If
\begin{equation}
\biggl(H-\frac{\lambda}{2}\biggl)^2\leq n+\frac{\lambda^2}{4},
\end{equation}
then $(H-\frac{\lambda}{2})^2\equiv n+\frac{\lambda^2}{4}$ and $M$ is isometric to the sphere $S^n(r)$ with radius $r>0$.
\end{proposition}

\section{Proofs of theorems for $\lambda$-hypersurfaces}
\noindent
In order to prove our theorems, we prepare several fundamental formulas.
Let $X: M^n\rightarrow\mathbb{R}^{n+1}$ be an
$n$-dimensional connected hypersurface of the $(n+1)$-dimensional Euclidean space
$\mathbb{R}^{n+1}$. We choose a local orthonormal frame field
$\{e_A\}_{A=1}^{n+1}$ in $\mathbb{R}^{n+1}$ with dual coframe field
$\{\omega_A\}_{A=1}^{n+1}$, such that, restricted to $M^n$,
$e_1,\cdots, e_n$ are tangent to $M^n$.
Then we have
\begin{equation*}
dX=\sum_i\limits \omega_i e_i, \quad de_i=\sum_j\limits \omega_{ij}e_j+\omega_{i n+1}e_{n+1}
\end{equation*}
and
\begin{equation*}
de_{n+1}=\sum_i\limits\omega_{n+1 i}e_i.
\end{equation*}
We restrict these forms to $M^n$, then
$$
\omega_{n+1}=0,\ \ \omega_{n+1i}=-\sum_{j=1}^nh_{ij}\omega_j,\ \ h_{ij}=h_{ji},
$$
where $h_{ij}$ denotes  components of the second fundamental form of $X: M^n\rightarrow\mathbb{R}^{n+1}$.
$H=\sum_{j=1}^nh_{jj}$
is the mean curvature and $II=\sum_{i,j}h_{ij}\omega_i\otimes\omega_jN$ is the second fundamental form
of $X: M^n\rightarrow\mathbb{R}^{n+1}$ with $N=e_{n+1}$.
Let
$$
h_{ijk}=\nabla_kh_{ij} \ \ {\rm and}\ \ h_{ijkl}=\nabla_l\nabla_kh_{ij},
$$
where $\nabla_j $ is the covariant differentiation operator.
Gauss equations,  Codazzi equations and Ricci formulas  are given by
\begin{equation}\label{eq:12-6-4}
R_{ijkl}=h_{ik}h_{jl}-h_{il}h_{jk},
\end{equation}
\begin{equation}\label{eq:12-6-5}
h_{ijk}=h_{ikj},
\end{equation}
\begin{equation}\label{eq:12-6-6}
h_{ijkl}-h_{ijlk}=\sum_{m=1}^nh_{im}R_{mjkl}+\sum_{m=1}^nh_{mj}R_{mikl},
\end{equation}
where $R_{ijkl}$ is  components of the curvature tensor.
For a function $F$, we denote covariant derivatives of $F$ by
$F_{,i}=\nabla_iF, \  F_{,ij}=\nabla_j\nabla_iF$.
For  $\lambda$-hypersurfaces,  an elliptic operator $\mathcal{L}$ is given  by
\begin{equation}\label{eq:12-6-1}
\mathcal{L}f=\Delta f-\langle X,\nabla f\rangle,
\end{equation}
where $\Delta$ and $\nabla$ denote the Laplacian and the gradient operator of the $\lambda$-hypersurface, respectively.
The $\mathcal{L}$ operator is introduced by Colding and Minicozzi in \cite{[CM]} for self-shrinkers and by Cheng and Wei \cite{cw2}
for $\lambda$-hypersurfaces.

\noindent
The following lemma\label{lemma 1} due to Colding and Minicozzi \cite{[CM]} is needed in order to prove our results.
\begin{lemma}\label{lemma 1}
Let $X: M\rightarrow \mathbb{R}^{n+1}$ be a complete hypersurface. If $u$, $v$ are $C^2$ functions satisfying
\begin{equation}
\int_M(|u\nabla v|+|\nabla u||\nabla v|+|u\mathcal{L}v|)e^{-\frac{|X|^2}{2}}d\mu< +\infty,
\end{equation}
then we have
\begin{equation}
\int_M u(\mathcal{L}v)e^{-\frac{|X|^2}{2}}d\mu=-\int_M\langle \nabla u,\nabla v\rangle e^{-\frac{|X|^2}{2}}d\mu.
\end{equation}
\end{lemma}

\vskip 4mm
\noindent {\it Proof of theorem 1.1}.
Since $\langle X,N\rangle +H=\lambda$, one has
\begin{equation}\label{eq:224-3}
H_{,i}=\sum_j h_{ij}\langle X,e_j\rangle,
\end{equation}
\begin{equation*}
H_{,ik}=\sum_jh_{ijk}\langle X,e_j\rangle +h_{ik}+\sum_jh_{ij}h_{jk}(\lambda-H).
\end{equation*}
From the Codazzi equation (\ref{eq:12-6-5}), we infer
\begin{equation*}
\Delta H=\sum_iH_{,ii}=\sum_iH_{,i}\langle X,e_i\rangle +H+S(\lambda-H).
\end{equation*}
Hence, we get
\begin{equation}\label{eq:224-11}
\aligned
\mathcal{L}H&=\Delta H-\sum_i\langle X,e_i\rangle H_{,i}=H+S(\lambda-H),
\endaligned
\end{equation}
\begin{equation}\label{eq:224-12}
\aligned
\dfrac12 \mathcal{L}H^2&=|\nabla H|^2+H^2+S(\lambda-H)H.
\endaligned
\end{equation}
By making use of the Ricci formulas and the Gauss equations and the Codazzi equations,   we have
\begin{equation*}
\aligned
\mathcal{L}h_{ij}&=\Delta h_{ij}-\sum_k\langle X,e_k\rangle h_{ijk}\\
&=\sum_kh_{ijkk}-\sum_k\langle X,e_k\rangle h_{ijk}\\
&=(1-S)h_{ij}+\lambda\sum_kh_{ik}h_{kj}.
\endaligned
\end{equation*}
Therefore, we obtain
\begin{equation*}
\aligned
\frac{1}{2}\mathcal{L}S&=\frac{1}{2}\bigl\{\Delta \sum_{i,j}(h_{ij})^2-\sum_k\langle X,e_k\rangle \bigl(\sum_{i,j}(h_{ij}^2)\bigl)_{,k}\bigl\}\\
&=\sum_{i,j,k}h_{ijk}^2+(1-S)\sum_{i,j}h_{ij}^2+\lambda\sum_{i,j,k}h_{ik}h_{kj}h_{ji}\\
&=\sum_{i,j,k}h_{ijk}^2+(1-S)S+\lambda f_3,\\
\endaligned
\end{equation*}
where $f_3=\sum_{i,j,k}h_{ij}h_{jk}h_{ki}$.

\noindent
Taking $\{e_1, e_2, \cdots, e_n\}$ such that $h_{ij}=\lambda_i\delta_{ij}$ at a point $p$ and putting $\mu_i=\lambda_i-\frac Hn$,
we have
$$
f_3=\sum_i\lambda_i^3=\sum_i(\mu_i+\frac Hn)^3
=B_3+\frac3nHB+\frac1{n^2}H^3
$$
with $B=\sum_i\mu_i^2=S-\frac{H^2}n$ and $B_3=\sum_i\mu_i^3$.
Thus, we have
\begin{equation*}
\aligned
\frac{1}{2}\mathcal{L}B&=\frac{1}{2}\mathcal{L}S-\dfrac12 \mathcal{L}\frac{H^2}n\\
&=\sum_{i,j,k}h_{ijk}^2-\frac1n|\nabla H|^2+(1-S)S+\lambda f_3-\frac{H^2}n-S(\lambda-H)\frac{H}n\\
&=\sum_{i,j,k}h_{ijk}^2-\frac1n|\nabla H|^2+(1-B)B-\frac1nH^2B+\lambda B_3+\frac2n\lambda HB.\\
\endaligned
\end{equation*}
Since
$$
\sum_i\mu_i=0, \ \ \ \sum_i\mu_i^2=B,
$$
it is not hard to prove
$$
|B_3|\leq \dfrac{n-2}{\sqrt{n(n-1)}}B^{\frac32}
$$
and the equality  holds if and only if at least,  $n-1$ of $\mu_i$  are equal.
Thus, we have
\begin{equation*}
\aligned
\frac{1}{2}\mathcal{L}B
&\geq \sum_{i,j,k}h_{ijk}^2-\frac1n|\nabla H|^2\\
&+(1-B)B-\frac1nH^2B-|\lambda| \dfrac{n-2}{\sqrt{n(n-1)}}B^{\frac32}+\frac2n\lambda HB\\
&=\sum_{i,j,k}h_{ijk}^2-\frac1n|\nabla H|^2\\
&+B\biggl(1-B-\frac1nH^2-|\lambda| \dfrac{n-2}{\sqrt{n(n-1)}}B^{\frac12}+\frac2n\lambda H\biggl)\\
&=\sum_{i,j,k}h_{ijk}^2-\frac1n|\nabla H|^2\\
&+B\biggl(1+\dfrac{n\lambda^2}{4(n-1)}-\frac1n(H-\lambda)^2-(\sqrt B+|\lambda| \dfrac{n-2}{2\sqrt{n(n-1)}})^2\biggl).\\
\endaligned
\end{equation*}
Since  $X: M\rightarrow \mathbb{R}^{n+1}$ has polynomial area growth, according to the results  of Cheng and Wei in \cite{cw2},
we can apply the lemma 2.1 to functions $1$ and $B=S-\frac{H^2}n$. Hence, we   have
\begin{equation*}
\aligned
0&\geq \int_{M}\biggl\{\sum_{i,j,k}h_{ijk}^2-\frac1n|\nabla H|^2\biggl\}e^{-\frac{|X|^2}{2}}d\mu\\
&+\int_{M}B\biggl(1+\dfrac{n\lambda^2}{4(n-1)}-\frac1n(H-\lambda)^2-(\sqrt B+|\lambda| \dfrac{n-2}{2\sqrt{n(n-1)}})^2\biggl)e^{-\frac{|X|^2}{2}}d\mu.
\endaligned
\end{equation*}
From the Codazzi equations and the Schwarz inequality, we have
$$
\sum_{i,j,k}h_{ijk}^2=3\sum_{i\neq k}h_{iik}^2+\sum_{i}h_{iii}^2+\sum_{i\neq j\neq k\neq i}h_{ijk}^2, \  \ \frac1n|\nabla H|^2\leq \sum_{i,k}h_{iik}^2,
$$
$$
\sum_{i,j,k}h_{ijk}^2-\frac1n|\nabla H|^2\geq 2\sum_{i\neq k}h_{ii k}^2+\sum_{i\neq j\neq k\neq i}h_{ijk}^2\geq 0
$$
and the equality holds if and only if  $h_{ijk}=0$ for any $i, j, k$.
Therefore, we get\\
\noindent either $B\equiv 0$ and  $X: M\rightarrow \mathbb{R}^{n+1}$ is totally umbilical;\\
\noindent  or
 there exists $p\in M$ such that
\begin{equation}
\biggl(\sqrt {S(p)-\dfrac{H^2(p)}n}+|\lambda| \dfrac{n-2}{2\sqrt{n(n-1)}}\biggl)^2+ \frac1n(H(p)-\lambda)^2> 1+\dfrac{n\lambda^2}{4(n-1)};
\end{equation}
or for any point of $M$
$$
\sum_{i,j,k}h_{ijk}^2-\frac1n|\nabla H|^2=0,
$$
$$
\biggl(\sqrt {S-\dfrac{H^2}n}-|\lambda| \dfrac{n-2}{2\sqrt{n(n-1)}}\biggl)^2+\frac1n(H-\lambda)^2 =1+\dfrac{n\lambda^2}{4(n-1)}.
$$
Hence,
we know that the second fundamental form is parallel, $X: M\rightarrow \mathbb{R}^{n+1}$ is  an isoparametric
complete hypersurface. If $\lambda=0$, then $X:M\rightarrow \mathbb{R}^{n+1}$ is isometric to the  sphere $S^n(\sqrt{n})$, the Euclidean space  $ \mathbb{R}^{n}$, the cylinder $S^k(\sqrt{k})\times   \mathbb{R}^{n-k}$.
If $\lambda\neq0$, then the number of the distinct principal curvatures is two and
one of them is simple,  $X:M\rightarrow \mathbb{R}^{n+1}$ is isometric to the  sphere $S^n(r)$,
 the Euclidean space  $ \mathbb{R}^{n}$, the cylinder $S^1(r)\times   \mathbb{R}^{n-1}$, the cylinder $S^{n-1}(r)\times   \mathbb{R}$.
The proof of theorem 1.1 is completed.
$$\eqno{\Box}$$

\vskip5mm
\noindent
By making use of the same assertions as in Cheng and Peng \cite{cp}, we know that
the following
generalized maximum principle holds.
\begin{theorem} {\rm(}Generalized maximum principle for $\mathcal{L}$-operator {\rm)}
Let $X: M^n\to \mathbb{R}^{n+1}$ be a complete $\lambda$-hypersurface  with Ricci
curvature bounded from below. Let $f$ be any $C^2$-function bounded
from above on this $\lambda$-hypersurface. Then, there exists a sequence of points
$\{p_k\}\subset M$, such that
\begin{equation}
\lim_{k\rightarrow\infty} f(X(p_k))=sup f,\quad
\lim_{k\rightarrow\infty} |\nabla f|(X(p_k))=0,\quad
\limsup_{k\rightarrow\infty}\mathcal{L} f(X(p_k))\leq 0.
\end{equation}
\end{theorem}

\vskip3mm
\noindent
{\it Proof of Theorem 1.3}.
From the proof in the theorem 1.1, we have
\begin{equation*}
\aligned
\frac{1}{2}\mathcal{L}B&\geq \sum_{i,j,k}h_{ijk}^2-\frac1n|\nabla H|^2\\
+&B\biggl(1+\dfrac{n\lambda^2}{4(n-1)}-\frac1n(H-\lambda)^2-(\sqrt B+|\lambda| \dfrac{n-2}{2\sqrt{n(n-1)}})^2\biggl)\\
\endaligned
\end{equation*}
and
$$
\sum_{i,j,k}h_{ijk}^2-\frac1n|\nabla H|^2\geq 2\sum_{i\neq k}h_{ii k}^2+\sum_{i\neq j\neq k\neq i}h_{ijk}^2\geq 0.
$$
Hence, we obtain
\begin{equation*}
\aligned
&\frac{1}{2}\mathcal{L}B\geq
B\biggl(1+\dfrac{n\lambda^2}{4(n-1)}-\frac1n(H-\lambda)^2-(\sqrt B+|\lambda| \dfrac{n-2}{2\sqrt{n(n-1)}})^2\biggl).
\endaligned
\end{equation*}
Since
$$
\sup \biggl\{(\sqrt {S-\dfrac{H^2}n}+|\lambda| \dfrac{n-2}{2\sqrt{n(n-1)}})^2+\frac1n(H-\lambda)^2 \biggl\}<1+\dfrac{n\lambda^2}{4(n-1)},
$$
we know $H^2$ and $S$ are bounded. Hence, from the Gauss equations, we infer
that the Ricci curvature is bounded from below. Applying  the generalized maximum principle for $\mathcal L$ of $\lambda$-hypersurfaces
to function $B$, there exists a sequence of points $\{p_k\} \subset  M$ such that
\begin{equation*}
\aligned
0\geq
\sup B\biggl(1+\dfrac{n\lambda^2}{4(n-1)}-\sup \biggl\{\frac1n(H-\lambda)^2+(\sqrt B+|\lambda| \dfrac{n-2}{2\sqrt{n(n-1)}})^2\biggl\}\biggl).
\endaligned
\end{equation*}
Hence, $\sup B=0$, that is, $S\equiv \frac{H^2}n$ and
$X: M\rightarrow \mathbb{R}^{n+1}$ is isometric to
\begin{enumerate}
\item the sphere $S^n(r)$ with radius $0<r<\sqrt{n}$ or
\item the Euclidean space  $ \mathbb{R}^{n}$.
\end{enumerate}
It completes the proof of the theorem 1.3.
$$
\eqno{\Box}
$$

\vskip3mm
\noindent
{\it Proof of Theorem 1.2}.
By a direct calculation, one obtains
\begin{equation}\label{eq:002}
\aligned
\frac{1}{2}\Delta |X|^2&=<\Delta X,X>+\sum_i<X_{,i},X_{,i}>\\
&=H<N,X>+n\\
&=n+\frac{\lambda^2}{4}-(H-\frac{\lambda}{2})^2.
\endaligned
\end{equation}
Since the assumption of polynomial area growth, we have
$$
\int_M(\Delta |X|^2)e^{-\frac{|X|^2}{2}}d\mu<+\infty,\ \ \int_M|\nabla |X|^2|^2e^{-\frac{|X|^2}{2}}d\mu<+\infty,
$$
then we can apply the lemma 2.1 to function $1$ and $|X|^2$ and obtain
\begin{equation*}
\frac{1}{4}\int_M|\nabla |X|^2|^2e^{-\frac{|X|^2}{2}}d\mu=\frac{1}{2}\int_M(\Delta |X|^2)e^{-\frac{|X|^2}{2}}d\mu
=\int_M(n+\frac{\lambda^2}{4}-(H-\frac{\lambda}{2})^2)e^{-\frac{|X|^2}{2}}d\mu.
\end{equation*}
From $(H-\frac{\lambda}{2})^2\geq n+\frac{\lambda^2}{4}$, we get
\begin{equation}
\biggl(H-\frac{\lambda}{2}\biggl)^2= n+\frac{\lambda^2}{4},\ \ \ \ \ \ <X,X>=r^2,
\end{equation}
namely,  $M$ is isometric to the sphere $S^n(r)$ with radius $r>0$.
It completes the proof of the proposition 1.2.
$$
\eqno{\Box}
$$

\vskip3mm
\noindent
{\it Proof of Proposition 1.1}.
Integrating \eqref{eq:002} over $M$ and using the Stokes formula, one concludes
\begin{equation}
\int_M (n+\frac{\lambda^2}{4}-(H-\frac{\lambda}{2})^2)d\mu=0,
\end{equation}
then it follows from $(H-\frac{\lambda}{2})^2\leq n+\frac{\lambda^2}{4}$ that
\begin{equation}
(H-\frac{\lambda}{2})^2= n+\frac{\lambda^2}{4}
\end{equation}
and $M$ is isometric to the sphere $S^n(r)$ with radius $r>0$.
It completes the proof of the proposition 1.1.
$$
\eqno{\Box}
$$

\end {document}